\documentclass[11pt]{amsart}
\usepackage{amssymb}
\usepackage{amsmath}
\usepackage{amsfonts}
\usepackage{amsthm}
\usepackage{amssymb}
\usepackage{mathrsfs}
\newtheorem{Theorem}{Theorem}[section]
\newtheorem{Proposition}[Theorem]{Proposition}

\newtheorem{Conjecture}[Theorem]{Conjecture}
\newtheorem{Corollary}[Theorem]{Corollary}
\theoremstyle{definition}
\newtheorem{Definition}[Theorem]{Definition}
\theoremstyle{remark}
\newtheorem{Remark}[Theorem]{Remark}

\DeclareMathOperator{\Area}{Area}

\newcommand{\dz}{\partial}
\newcommand{\dzb}{\bar\partial}
\newcommand{\ind}{\mathrm{Ind}}
\newcommand{\nul}{\mathrm{Nul}}
\begin{document}
\title[An inequality for Laplace eigenvalues on the
sphere]{An isoperimetric inequality 
for Laplace eigenvalues on the sphere}
\author[M. Karpukhin]{Mikhail Karpukhin\textsuperscript{1}}
\address[Mikhail Karpukhin]{Department of Mathematics and Statistics,
McGill University, Burnside Hall, 805 Sherbrooke Street West, Montr\'eal, 
Qu\'ebec, H3A 0B9, Canada}
\curraddr{Department of Mathematics, University of California, Irvine, 
CA 92697, USA}
\email{mkarpukh@uci.edu}
\author[N. Nadirashvili]{Nikolai Nadirashvili}
\address[Nikolai Nadirashvili]{CNRS, I2M UMR 7353 --- 
Centre de Math\'e\-ma\-ti\-ques et Informatique,
Marseille, France}
\email{nikolay.nadirashvili@univ-amu.fr}
\author[A. Penskoi]{Alexei V. Penskoi\textsuperscript{2}}
\address[Alexei V. Penskoi]{Department of Higher Geometry and Topology, 
Faculty of Mathematics and Mechanics, Moscow State University,
Leninskie Gory, GSP-1, 119991, Moscow, Russia \newline {\em and}
\newline Faculty of Mathematics,
National Research University Higher School of Economics,
6 Usacheva Str., 119048, Moscow, Russia \newline {\em and}
\newline Interdisciplinary Scientific Center
J.-V. Poncelet (ISCP, UMI 2615), Bolshoy Vlasyevskiy 
Pereulok 11, 119002, Moscow, Russia \newline {\em and}
\newline (Visiting position) CNRS, I2M UMR 7353 --- 
Centre de Math\'e\-ma\-ti\-ques et Informatique,
Marseille, France}
\email[Corresponding author]{penskoi@mccme.ru}
\author[I. Polterovich]{Iosif Polterovich\textsuperscript{3}}
\address[Iosif Polterovich]{D\'epartement de math\'ematiques et de 
statistique, Universit\'e de Montr\'eal, CP 6128 Succ. Centre-Ville, 
Montr\'eal, Qu\'ebec, H3C 3J7, Canada}
\email{iossif@dms.umontreal.ca}
\date{}

\footnotetext[1]{Partially supported
by Tomlinson Fellowship and Schulich Fellowship.}
\footnotetext[2]{Partially supported by the Simons-IUM Fellowship and 
by the Young Russian Mathematics award.}
\footnotetext[3]{Partially supported by 
NSERC, FRQNT and the Canada Research Chairs 
program.}
\begin{abstract}
We show that for any positive integer $k$, the $k$-th 
nonzero eigenvalue of the Laplace-Beltrami 
operator on the two-dimensional sphere endowed with a 
Riemannian metric of unit area, is maximized 
in the limit by a sequence of metrics converging to 
a union of $k$ touching identical round spheres. 
This proves a conjecture posed by the second author in 
2002 and yields a sharp isoperimetric 
inequality for all nonzero eigenvalues of the Laplacian on a 
sphere. Earlier, the result was known 
only for $k=1$ (J.~Hersch, 1970), $k=2$ (N.~Nadirashvili, 2002; 
R.~Petrides, 2014) and $k=3$ 
(N.~Nadirashvili and Y.~Sire, 2017). In particular, we argue  that for any 
$k\geqslant 2$,  the supremum of the $k$-th nonzero eigenvalue 
on a sphere of unit area
is not attained in the class of  Riemannian metrics which are smooth 
outsitde a finite set of  conical singularities.
The proof  uses certain properties of harmonic  maps between spheres, 
the key new ingredient being  a bound on the 
harmonic degree of a harmonic map into a sphere obtained by N. Ejiri.
\end{abstract}
\maketitle
\section{Introduction and main results}
Let $M$ be a closed surface and $g$ be
a Riemannian metric on $M.$
Consider the Laplace-Beltrami
operator $\Delta:C^\infty(M)\longrightarrow C^\infty(M)$
associated with the metric $g,$
$$
\Delta f=-\frac{1}{\sqrt{|g|}}\frac{\partial}{\partial x^i}%
\left(\sqrt{|g|}g^{ij}\frac{\partial f}{\partial x^j}\right),
$$
and its eigenvalues
\begin{equation}\label{eigenvalues}
0=\lambda_0(M,g)<\lambda_1(M,g)\leqslant%
\lambda_2(M,g)\leqslant\lambda_3(M,g)\leqslant\dots,
\end{equation}
where each eigenvalue is repeated according to its multiplicity.

Let $\Area(M,g)$ denote the area of $M$ with respect to the
Riemannian metric $g.$
Since the eigenvalues possess a rescaling property,
$$
\forall t>0\quad\lambda_i(M,tg)=\frac{\lambda_i(M,g)}{t},
$$
it is natural to consider functionals
$$
\bar{\lambda}_i(M,g)=\lambda_i(M,g)\Area(M,g),
$$
invariant under the rescaling transformation $g\mapsto tg.$
The functional $\bar{\lambda}_i(M,g)$
is sometimes called an eigenvalue normalized by the area or 
simply a normalized eigenvalue.

Our principal interest is in the  following geometric optimization problem 
for the Laplace-Beltrami eigenvalues:
given a compact surface $M$ and a number $i\in\mathbb{N}$
find the quantity
\begin{equation}\label{Lambda-definition}
\Lambda_i(M)=\sup_g\bar\lambda_i(M,g),
\end{equation}
where the supremum is taken over the space of all
Riemannian metrics $g$ on $M.$

Since the functional $\bar{\lambda}_i(M,g)$
is invariant under rescaling of the metric 
$g\mapsto tg,$ where $t\in\mathbb{R}_+,$ this problem is
equivalent to the question of finding $\sup\lambda_i(M,g),$
where the supremum is taken over the space of all
Riemannian metrics $g$ of area $1$ on $M.$
This problem is  often referred to as  the {\em isoperimetric
problem for the Laplace-Beltrami eigenvalues}, since  
for any metric $g$ of unit area on $M$,
one has the inequality
$$
\lambda_i(M,g)\leqslant\Lambda_i(M).
$$

\begin{Definition}\label{maximal-closed}
Let $M$ be a closed surface.
A metric $g_0$ on $M$ is called {\it maximal} for the functional
$\bar{\lambda}_i(M,g)$ if
$$
\Lambda_i(M)=\bar{\lambda}_i(M,g_0).
$$
\end{Definition}

Note that if a maximal metric exists, it is defined up to
multiplication by a positive constant due to the rescaling
invariance of the functional.

The main goal of the present paper is to provide a complete solution of the 
isoperimetric problem for the Laplace-Beltrami eigenvalues on a sphere.
We will say that a metric on a sphere is 
{\it round}\/ if it is isometric to the standard round metric
of constant Gau\ss{}ian curvature. We will also say that 
the two round spheres are {\it touching} if they have exactly 
one point in common. In what follows, the choice of  the points 
where the touching occurs is of no significance, as it does not 
change the  spectrum (note the spectra of touching and disjoint 
spheres are the same). 
\begin{Theorem}\label{maintheorem} Let  $(\mathbb{S}^2, g)$ be 
the sphere endowed with a  metric $g$  which is smooth outside 
a finite number of conical singularities. Then the following  inequality holds: 
\begin{equation}
\label{main:isop}
\bar \lambda_k(\mathbb{S}^2,g)\leqslant 8\pi k, \,\,\, k\geqslant 1.
\end{equation}
For $k=1$,  the equality is attained if and only if $g$ is a
round metric. For $k\geqslant 2$ the inequality is strict, and the equality 
is attained in the limit by a  sequence of metrics degenerating to a 
union of $k$ touching identical round spheres.
\end{Theorem}
Theorem~\ref{maintheorem} immediately implies 
\begin{equation}
\label{main}
\Lambda_k(\mathbb{S}^2)=8\pi k, \,\,\, k\geqslant 1.
\end{equation}
It settles a conjecture that was proposed in 2002 by the second 
author~\cite{Nadirashvili2002} (see also~\cite{Nadirashvili-Sire2015q2}, 
as well as \cite{Braxton} for numerical evidence supporting the conjecture). 
For $k=1$, the result 
is a classical theorem  due to J.~Hersch~\cite{Hersch1970}. For $k=2$,  the 
proof of the 
equality~\eqref{main}  was outlined in~\cite{Nadirashvili2002}, and a 
different argument was 
found later in~\cite{Petrides2014}. Recently,~\eqref{main} has been 
also established in the case 
$k=3$~\cite{Nadirashvili-Sire2015q2}.  For  $k\geqslant 4$, 
Theorem~\ref{maintheorem} is new.   
\begin{Remark}
The degenerating sequence of metrics  in Theorem \ref{maintheorem} 
illustrates the "bubbling phenomenon'' arising in the maximization 
of higher eigenvalues, see \cite{Nadirashvili-Sire2015b} for details. 
It would be interesting to check whether the equality in \eqref{main:isop} 
could be attained in the limit only by a sequence of metrics converging to  
a union of touching spheres,
or there exist other maximizing sequences.
\end{Remark}
\begin{Remark} It is instructive to compare Theorem~\ref{maintheorem} 
with isoperimetric inequalities for other eigenvalue problems. In particular, 
one could draw an analogy with the Hersch-Payne-Schiffer inequality for 
the $k$-th nonzero Steklov eigenvalue on simply connected planar domains, 
which is sharp for all $k\geqslant 1$, and the equality is attained in the limit 
as the domain degenerates into a disjoint union of $k$ identical 
disks~\cite{GP}. The situation is quite different for higher  Dirichlet 
and Neumann  eigenvalues on planar domains. In this case, sharp 
isoperimetric inequalities for eigenvalues for $k\geqslant 3$ are not known,  
and numerical experiments indicate  that optimal shapes  may have rather complicated  geometries \cite{AF}. Moreover, the celebrated Polya's 
conjecture (see \cite{Pol1, Pol2, LY1983}) suggests that the  $k$-th 
Dirichlet (respectively, nonzero  Neumann) eigenvalue is bounded 
below (respectively, above) for all $k\geqslant 1$ by the main term in 
Weyl's asymptotics,   which is not the case for inequality \eqref{main:isop} 
(see \cite[Remark 1.2.8]{GNP} for a further discussion).
\end{Remark}

\subsection*{Organization of the paper}
In the next section we give a brief overview of results on 
isoperimetric inequalities and extremal metrics for 
Laplace-Beltrami eigenvalues on surfaces.
In Section~\ref{minimal_immersions} we recall the relation
between extremal metrics, minimal immersions into spheres 
and harmonic maps. We also extend our considerations
to harmonic immersions with branch points and metrics with conical singularities.
In Section~\ref{Calabi-Barbosa-implications} we state the Calabi--Barbosa
theorem on harmonic immersions $\mathbb{S}^2\longrightarrow\mathbb{S}^n$ 
with branch points, and define the notion of the harmonic degree.
Up to this point, our exposition follows 
closely~\cite{Nadirashvili-Penskoi2016q} where similar ideas have been 
used to obtain a sharp isoperimetric inequality for the second eigenvalue on 
the projective plane. The main novel ingredient of our approach is 
introduced at the end of Section~\ref{Calabi-Barbosa-implications}. 
It is an upper bound due to N. Ejiri~\cite{Ejiri1998} on the harmonic degree 
of a harmonic immersion $f:\mathbb{S}^2\longrightarrow\mathbb{S}^{2m}$
in terms of the index $k$ of the eigenvalue $\lambda_k$ extremized by the
metric induced by $f$.
The proof of
Theorem \ref{maintheorem} is presented in  Section~\ref{mainproof}. 
For the convenience of the reader,  in Appendix \ref{sec:Ejiri} we 
sketch the  proof of Ejiri's bound, which  uses techniques  from 
algebraic geometry. These methods are relatively new in  spectral 
geometry and appear to  be significant  in the study of 
extremal metrics for Laplace-Beltrami eigenvalues 
(see also~\cite{Nayatani-Shoda2017q}).

\smallskip

\noindent {\bf Acknowledgements.}
A.P. is very grateful to 
the Institut de Ma\-th\'e\-ma\-ti\-ques de Marseille (I2M, UMR~7373)
for its hospitality.
The authors would like to thank G.~Kokarev, E.~Lauret  and the 
anonymous referees for helpful remarks.

\section{Isoperimetric inequalities and extremal metrics 
for eigenvalues on surfaces}

In this section we present a brief overview of the main 
developments in the study of the 
isoperimetric inequalities and extremal metrics for the 
eigenvalues of the Laplace-Beltrami 
operator on surfaces (see also~\cite[section III.8]{Schoen-Yau1994}).
The story began almost half a century 
ago when Hersch~\cite{Hersch1970}  
obtained a sharp isoperimetric inequality for the first 
eigenvalue on a sphere which was 
mentioned in the previous section.  Later, Yang and Yau~\cite{Yang-Yau1980} 
proved that
if $M$ is an orientable surface of genus $\gamma$ then
$$
\bar{\lambda}_1(M,g)\leqslant 8\pi(\gamma+1).
$$
Actually, the arguments of Yang and Yau imply a stronger, 
though in general still non-sharp estimate, 
$$
\bar{\lambda}_1(M,g)\leqslant 8\pi\left[\frac{\gamma+3}{2}\right],
$$
see~\cite{ElSoufi1984}.
Here $[\cdot]$ denotes the integer part of a number.

A similar result in the non-orientable case has been recently 
obtained in~\cite{Karpukhin2016}, see also~\cite{Kokarev2017q}.
For a non-orientable surface $M$ of genus $\gamma,$ 
$$
\bar\lambda_1(M,g)\leqslant 16\pi\left[\frac{\gamma+3}{2}\right].
$$
Here the genus of a non-orientable surface is defined to be 
the genus of 
its orientable double covering.

For higher eigenvalues, Hassannezhad~\cite{Hassannezhad2013}, 
building on the earlier work of Korevaar ~\cite{Korevaar1993}, proved 
that there exists a constant $C$ such that for any 
$i>0$ and any compact surface $M$ of genus $\gamma$,
the following upper bound holds:
$$
\bar{\lambda}_i(M,g)\leqslant C(\gamma+i).
$$

These results imply that the functionals $\bar{\lambda}_i(M,g)$
are bounded from above, and, as a result, $\Lambda_i(M)<+\infty.$
A number of important results on the existence and regularity 
of maximizing metrics has been recently obtained, 
see~\cite{Kokarev2014a, Nadirashvili-Sire2015a, 
Nadirashvili-Sire2015b, Petrides2014b}. However, the list of known values  
$\Lambda_i(M)$ for all surfaces other than the sphere is 
quite short.

Li and Yau proved in 1982 in the paper~\cite{Li-Yau1982}
that the standard metric on the projective plane
is the unique maximal metric for
$\bar{\lambda}_1(\mathbb{RP}^2,g)$
and 
$$
\Lambda_1(\mathbb{RP}^2)=12\pi.
$$

The second author proved in 1996 in the paper~\cite{Nadirashvili1996}
that the standard metric on the equilateral torus
is the unique maximal metric for
$\bar{\lambda}_1(\mathbb{T}^2,g)$
and 
$$
\Lambda_1(\mathbb{T}^2)=\frac{8\pi^2}{\sqrt{3}}.
$$

Note that the functional $\bar{\lambda}_i(M,g)$ depends continuously
on the metric $g.$ 
However, when $\bar{\lambda}_i(M,g)$ is a multiple eigenvalue 
this functional is not in general differentiable. If we consider an
analytic variation $g_t$ of the metric $g=g_0,$
then it was proved by Berger~\cite{Berger1973},
Bando and Urakawa~\cite{Bando-Urakawa1983},
El Soufi and Ilias~\cite{ElSoufi-Ilias2008}
that the left and right derivatives
of the functional $\bar{\lambda}_i(M,g_t)$ with respect to $t$ exist.
This leads us to the  following definition
given by the first author in the paper~\cite{Nadirashvili1996}
and by El Soufi and Ilias in the 
papers~\cite{ElSoufi-Ilias2000,ElSoufi-Ilias2008}.

\begin{Definition}
A Riemannian metric $g$ on  a closed surface
$M$ is called an {\it extremal metric} for the
functional $\bar\lambda_i(M,g)$ (or simply for the eigenvalue $\lambda_i$) if for any analytic deformation
$g_t$ such that $g_0=g$ one has
$$
\frac{d}{dt}\bar{\lambda}_i(M,g_t)%
\left.\vphantom{\raisebox{-0.5em}{.}}\right|_{t=0+}\, %
\cdot\, \frac{d}{dt}\bar{\lambda}_i(M,g_t)%
\left.\vphantom{\raisebox{-0.5em}{.}}\right|_{t=0-} \leqslant 0.
$$
\end{Definition}

It was proved in~\cite{Jakobson-Nadirashvili-Polterovich2006}
that the metric on the Klein bottle realized as the so-called
bipolar Lawson surface $\tilde{\tau}_{3,1},$
is extremal for $\bar\lambda_1(\mathbb{KL},g).$
It was conjectured in this paper that this metric 
is the unique extremal metric and, moreover,
the maximal one. 
It was later shown by El Soufi, Giacomini and 
Jazar~\cite{ElSoufi-Giacomini-Jazar2006}
that this metric on $\tilde{\tau}_{3,1}$
is indeed the unique extremal metric for the functional
$\bar\lambda_1(\mathbb{KL},g).$
Moreover, one can prove that $\Lambda_1(\mathbb{KL})$ is attained on 
a smooth metric~\cite{Cianci-Medvedev}, and therefore the metric 
on $\tilde{\tau}_{3,1}$ is the maximal one.
Hence,
$$
\Lambda_1(\mathbb{KL})=\bar{\lambda}_1(\mathbb{KL},g_{\tilde{\tau}_{3,1}})=%
12\pi E\left(\frac{2\sqrt{2}}{3}\right),
$$
where $E$ is the complete elliptic integral of the second kind
and $g_{\tilde{\tau}_{3,1}}$ is the metric on $\tilde{\tau}_{3,1}.$

More results on extremal metrics on tori and Klein bottles
could be found in the 
papers~\cite{ElSoufi-Ilias2000,Karpukhin2013,Karpukhin2014,
Karpukhin2015,Lapointe2008,Penskoi2012,Penskoi2013,Penskoi2015}. 
A review of these results is given by the third author in 
the paper~\cite{Penskoi2013a}.

It was shown in \cite{JLNNP} using a combination of analytic and numerical tools
that the maximal metric for the first eigenvalue on 
the surface of genus two is the metric on the Bolza surface $\mathcal P$
induced from the canonical metric on the sphere
using the standard covering ${\mathcal P}\longrightarrow\mathbb{S}^2.$ 
The authors stated this result as a conjecture, because
the argument is partly  based on a numerical calculation. The proof of
this conjecture was presented  in a recent preprint ~\cite{Nayatani-Shoda2017q}. 

For surfaces of  genus $\gamma\geqslant 3$, no maximal metrics for 
Laplace eigenvalues are known.
The only surface other than the sphere for which a sharp 
isoperimetric inequality was 
established for a higher eigenvalue is the projective plane.
It was shown in~\cite{Nadirashvili-Penskoi2016q} that
$$
\Lambda_2(\mathbb{RP}^2)=20\pi.
$$
It turns out that there is no maximal metric,
but the supremum could be
obtained as a limit on a sequence of metrics converging to a union 
of a projective plane and a sphere touching
at a point, each endowed with the standard metric,  
such that the ratio of the areas 
of the projective plane 
and the sphere is $3:2.$ Inspired by this result, as 
well as by Theorem~\ref{maintheorem}, we 
propose  the following conjecture.
\begin{Conjecture}
The equality
$$
\Lambda_k(\mathbb{RP}^2)=4\pi (2k+1)
$$
holds for any $k\geqslant 1$. For $k\geqslant 2$ the supremum can not be 
attained on a smooth metric, and is realized in the limit by a 
sequence of metrics degenerating  to a union of $k-1$ 
identical round spheres and a standard 
projective plane  touching each other,
such that the ratio of the areas of the projective plane 
and the spheres is $3:2$. 
\end{Conjecture}

Let us conclude this section by a result  which  is a special 
case  of \cite[Theorem 2]{Petrides2017};  it could be also 
deduced from  \cite[Theorem 1.1]{Nadirashvili-Sire2015b}. 
We use here that there exists only one conformal structure
on $\mathbb{S}^2$.
\begin{Theorem}
\label{thm:petrides} Let $k\geqslant 2$ and suppose that the 
following strict inequality holds:
\begin{equation}
\label{petrides}
\Lambda_k(\mathbb{S}^2) > \max_{\substack{1<s\leqslant k\\%
k_1+\dots +k_s=k\\k_1,\dots,k_s\geqslant 1}}\,\,  \sum_{j=1}^s %
\Lambda_{k_j}(\mathbb{S}^2).
\end{equation}
Then there exists a maximal metric $g$, which is smooth except 
possibly at a finite number of conical singularities, such that
$\Lambda_k(\mathbb{S}^2)=\bar\lambda_k(\mathbb{S}^2,g)$.
\end{Theorem}
Theorem \ref{thm:petrides} will be used in an essential  
way in the  proof of Theorem~\ref{maintheorem}, see Section \ref{mainproof}.
\begin{Remark}
\label{gluing}
As mentioned in \cite{Petrides2017},  if one replaces the strict 
inequality in \eqref{petrides}  by a non-strict inequality, then 
by a gluing argument it is always true. Indeed, 
assume that the maximum $\mu_k$ of the right hand 
side is achieved for some $k_1,\dots, k_s$. Take $s$ copies of the 
sphere, and on the $j$-th copy, $j=1,\dots, s$,  consider a metric 
$g_j$ such that $\bar\lambda_{k_j}(\mathbb{S}^2,g_j)$ is close 
to $\Lambda_{k_j}(\mathbb{S}^2)$. Scale $g_j$ by a positive constant 
$t_j$ such that $\lambda_{k_j}(\mathbb{S}^2, t_jg_j)$ is independent 
of $j$. Gluing these spheres together using sufficiently thin and short 
tubes  yields a metric $g$ on the sphere with $\bar\lambda_k(\mathbb{S}^2,g)$ 
close to $\mu_k$ (see \cite[Section 2]{CE} for details). Since $g_j$ could be 
chosen in such a way that  $\bar\lambda_{k_j}(\mathbb{S}^2,g_j)$ is 
arbitrary close to $\Lambda_{k_j}(\mathbb{S}^2)$ for all $j=1,\dots, s$, 
and the  connecting tubes could be made arbitrary small, we 
get  $\Lambda_k(\mathbb{S}^2)\geqslant \mu_k$ which proves the claim.
\end{Remark}

\begin{Remark}
In fact, both~\cite[Theorem 2]{Petrides2017} 
and~\cite[Theorem 1.1]{Nadirashvili-Sire2015b} are concerned with 
the existence of Riemannian metrics with conical singularities for 
which the quantity 
$$
\Lambda_k(M,[g]) = \sup_{h\in [g]}\bar\lambda_k(M,h)
$$
is attained.
The existence is guaranteed provided the 
condition~\cite[formula (0.3)]{Petrides2017} holds.
Using Theorem~\ref{maintheorem} one can simplify this condition, 
see discussion in~\cite[page 5]{Petrides2017}:
\begin{Corollary}
Let $(M,g)$ be a Riemannian surface and $k\geqslant 2$. If 
$$
\Lambda_k(M,[g])>\Lambda_{k-1}(M,[g]) + 8\pi,
$$
then there exists a maximal metric $\tilde g\in [g]$, smooth 
except possibly at a finite set of conical singularities, such 
that $\Lambda_k(M,[g]) = \bar\lambda_k(M,\tilde g)$. 
\end{Corollary}
\end{Remark}

\section{Extremal metrics, minimal immersions and harmonic maps}
\label{minimal_immersions}

Let us recall the definition of a minimal map, see 
e.g.~\cite{Eells-Lemaire1978,Eells-Sampson1964}.

\begin{Definition}
Let $M$ and $N$ be smooth manifolds and $h$ be a Riemannian metric on $N$.
An immersion $f:M\looparrowright N$ is called
{\it minimal} 
$f$ if it is extremal for the volume functional
$$
V[f]=\int_M\,dVol_{f^*h}.
$$
If $M$ is a Riemannian manifold endowed with a metric 
$g=f^*h,$ then $f$ is called {\em an isometric minimal immersion}
and $M$  is referred to as {\em an (immersed) minimal submanifold}.
\end{Definition}
It is well-known that a submanifold 
$M\looparrowright\mathbb{R}^n$ is minimal if and only if
the coordinate functions $x^i$ are harmonic
with respect to the Laplace-Beltrami operator on $M.$
Since harmonic functions are eigenfunctions
with eigenvalue~$0,$ it is natural to ask what is an analogue
of this statement for a non-zero eigenvalue. The answer was
given by Takahashi in 1966.

\begin{Theorem}[Takahashi~\cite{Takahashi1966}]\label{takahashi}
If an isometric immersion 
$$
f:M\looparrowright\mathbb{R}^{n+1},\quad f=(f^1,\dots,f^{n+1}),
$$
is defined by eigenfunctions
$f^i$ of the Laplace-Beltrami operator $\Delta$
with a common eigenvalue $\lambda,$
$$
\Delta f^i=\lambda f^i,
$$
then 

\noindent (i) the image $f(M)$ lies on the sphere $\mathbb{S}^n_R$
of radius $R$ with the center at the origin
such that 
\begin{equation}\label{lambda-R}
\lambda=\frac{\dim M}{R^2},
\end{equation}

\noindent (ii) the immersion $f:M\looparrowright\mathbb{S}^n_R$
is minimal.

If 
$$
f:M\looparrowright\mathbb{S}^{n}_R\subset\mathbb{R}^{n+1},
\quad f=(f^1,\dots,f^{n+1}),
$$ 
is an isometric minimal immersion of a manifold $M$
into the sphere $\mathbb{S}^{n}_R$
of radius $R,$ then $f^i$
are eigenfunctions of the Laplace-Beltrami
operator~$\Delta,$
$$
\Delta f^i=\lambda f^i,
$$
with the same eigenvalue $\lambda$
given by formula~\eqref{lambda-R}.
\end{Theorem}

Let us introduce the counting function for the 
Laplace eigenvalues of a manifold $(M,g)$:
\begin{equation}
\label{counting}
N(\lambda; M,g)=\#\{\lambda_i(M,g)<\lambda\}.
\end{equation}
Theorem~\ref{takahashi} implies that if $M$
is isometrically minimally immersed in the sphere $\mathbb{S}^{n}_R,$
then the eigenvalue $\frac{\dim M}{R^2}$ is multiple. It is easy
to see that $\lambda_{N\left(\frac{\dim M}{R^2}; M,g\right)} (M,g)$ is the first
Laplace eigenvalue equal to $\frac{\dim M}{R^2}$, where $g$ is the  
metric on $M$ induced by the immersion.  This fact  is important due to the
following theorem.

\begin{Theorem}[Nadirashvili~\cite{Nadirashvili1996},
El Soufi, Ilias~\cite{ElSoufi-Ilias2008}]\label{Nadirashvili-ElSoufi-Ilias}
Let $M\looparrowright\mathbb{S}^n_R$ be an immersed minimal
compact submanifold. Then the metric $g$ induced on $M$ by this immersion
is extremal for the eigenvalue
$\lambda_{N\left(\frac{\dim M}{R^2};M,g\right)}$.

If a metric $g$ on a compact manifold $M$ is extremal for 
the $k$-th eigenvalue, then $g$ is induced by
an isometric minimal immersion $M\looparrowright\mathbb{S}^n_R$  to the sphere
of some dimension $n\geqslant 2$ and radius $R=\sqrt{\dim M/\lambda_k(M,g)}$.
Moreover, the components of the immersion are given by 
the $\lambda_k$-eigenfunctions of the Laplacian.
\end{Theorem}

We later use  Theorem~\ref{Nadirashvili-ElSoufi-Ilias} for 
$M=\mathbb{S}^2.$ In this case $\dim M=2.$ Using rescaling one can consider only
the case $R=1.$ Hence, the value
$$
N\left(\frac{\dim\mathbb{S}^2}{R^2}; \mathbb{S}^2, g \right)=N(2;\, \mathbb{S}^2, g)
$$
is particularly significant for our further considerations.

Let us now  recall the definition of a harmonic map,
see e.g. the review~\cite{Eells-Lemaire1978}.

\begin{Definition}
Let $(M,g)$ and $(N,h)$ be Riemannian manifolds. A smooth
map $f:M\longrightarrow N$ is called
{\it harmonic} if $f$ is extremal for the energy functional
\begin{equation}\label{energy-functional}
E[f]=\int_M|df(x)|^2\,dVol_g.
\end{equation}
Here, $df$ is viewed as a section of $T^*M\otimes f^*TN$ with 
the metric induced by $g$ and $h$. In local coordinates
$$
|df|^2 = g^{kl}\frac{\partial f^i}{\partial x^k}%
\frac{\partial f^j}{\partial x^l}h_{ij}.
$$
\end{Definition}

The following theorem (see, e.g. the paper~\cite{Eells-Sampson1964})
explains the relation between minimal and harmonic maps
in the class of isometric immersions.

\begin{Theorem}\label{harmonic-minimal}
Let $M,$ $N$ be Riemannian manifolds. If $f:M\looparrowright N$
is an isometric immersion, then $f$ is harmonic if and only
if $f$ is minimal.
\end{Theorem}

It is well-known that if $\dim M=2$, then the property of $f$
to be harmonic depends only on the conformal class of the
metric $g,$ see e.g.~\cite{Eells-Lemaire1978}.
A map $f:(M,g)\longrightarrow (N,h)$ is called {\it conformal} if 
the metrics $g$ and $f^*h$ belong to the same conformal class. 
Together with Theorem~\ref{harmonic-minimal} 
this yields the following corollary.

\begin{Corollary}
Let $(M,g)$ be a Riemannian surface. A conformal immersion 
$f\colon (M,g) \looparrowright (N,h)$ is harmonic iff 
$f\colon M \looparrowright (N,h)$ is minimal.
\end{Corollary}
It is a well-known fact that there is only one conformal
class of Riemannian metrics on $\mathbb{S}^2.$ Hence, for $M=\mathbb{S}^2$
any harmonic immersion is conformal. Therefore, 
Theorem~\ref{Nadirashvili-ElSoufi-Ilias} implies:

\begin{Corollary}\label{metrics-harmonic-maps-s2}
The smooth extremal metrics for the eigenvalues of the Laplace-Beltrami operator
on the sphere $\mathbb{S}^2$ are exactly the
metrics induced on $\mathbb{S}^2$ by harmonic
immersions $f:\mathbb{S}^2\looparrowright\mathbb{S}^n$, 
where $\mathbb{S}^n$ is endowed with a round metric.
\end{Corollary}

However, as one can see from Theorem~\ref{thm:petrides} it 
is not sufficient to consider smooth metrics as one has to account 
for a possible presence of conical singularities. This is done by 
allowing minimal immersions to have branch points. 
According to~\cite{Gulliver-Osserman-Royden1973}, it is equivalent to 
replacing conformal harmonic maps with weakly conformal harmonic maps.
\begin{Definition}
A map $f\colon (M,g)\looparrowright(N,h)$ is called {\em weakly conformal} 
if $f^*h = \alpha g$, where $\alpha\geqslant 0$ is a non-negative smooth 
function with isolated zeroes, i.e. the differential of $f$ can be 
degenerate at isolated points of $M$.
\end{Definition}  
For a more detailed discussion on minimal immersions with branch 
points and metrics with conical singularities, 
see~\cite{Cianci-Medvedev,Nadirashvili-Penskoi2016q}. In particular, 
Theorem~\ref{Nadirashvili-ElSoufi-Ilias} can be generalized to 
this context, see \cite[Corollary 4.7]{Kokarev2014a}. In the special 
case  $M=\mathbb{S}^2$ we have:
\begin{Theorem}
\label{metrics-harmonic-maps-s2-branch}
Let $f: \mathbb{S}^2 \looparrowright\mathbb{S}^n$ be a harmonic map. 
Then the metric $g$ (possibly with conical singularities) 
on $\mathbb{S}^2$ induced by this map
is extremal for the eigenvalue
$\lambda_{N\left(2;\, \mathbb{S}^2,g\right)}$.

If a metric $g$ (possibly with conical singularities)  
on  $\mathbb{S}^2$ is extremal for the $k$-th eigenvalue,  
then $g$ coincides up to a rescaling with the metric  induced 
by a harmonic map  $\mathbb{S}^2 \looparrowright\mathbb{S}^n$ to the 
unit sphere of some dimension $n\geqslant2$. Moreover, 
$\lambda_k=2$ and the components of the harmonic map are given by 
the $\lambda_k$-eigenfunctions of the Laplacian.
\end{Theorem}

\begin{Remark} As was confirmed  to us by G. Kokarev,  
the ``if'' part of  \cite[Corollary 4.7]{Kokarev2014a} 
assumes additionally that $\lambda_{k-1}<\lambda_k$ or 
$\lambda_k<\lambda_{k+1}$, cf. \cite[Theorem 3.1]{ElSoufi-Ilias2008}. 
Note that by definition \eqref{counting} of the counting 
function, $\lambda_{N\left(2;\, \mathbb{S}^2,g\right)}$ satisfies 
the former inequality, i.e. it is strictly greater than the preceding eigenvalue.
\end{Remark}

Let us also remark that for any manifold equipped with a metric
with isolated conical singularities 
it is possible to construct
a sequence of smooth Riemannian manifolds such that 
their area as well as their eigenvalues
converge to the area and eigenvalues of the initial manifold,
see e.g.~\cite{Sher2015}.

\section{Harmonic maps into spheres and the harmonic degree}
\label{Calabi-Barbosa-implications}
We now focus on  harmonic immersions  $\mathbb{S}^2\looparrowright\mathbb{S}^n$
with branch points. 
Let us introduce the following useful definition.

\begin{Definition}
A harmonic map $f:\Sigma\longrightarrow\mathbb{S}^n\subset\mathbb{R}^{n+1}$
of a surface $\Sigma$ to the standard unit 
sphere $\mathbb{S}^n\subset\mathbb{R}^{n+1}$
is called {\it linearly full} if the image $f(\Sigma)$ does not lie in a hyperplane
of $\mathbb{R}^{n+1}.$
\end{Definition}

The following theorem was proved by Calabi in 1967 and later refined
by Barbosa in 1975. 
Let $g_{\mathbb{S}^n}$ denote the standard metric on $\mathbb{S}^n$ of radius $1.$

\begin{Theorem}[Calabi~\cite{Calabi1967}, 
Barbosa~\cite{Barbosa1975}]\label{Calabi-Barbosa}
Let $f:\mathbb{S}^2\longrightarrow\mathbb{S}^n$
be a linearly full harmonic immersion with branch points.
Then
\begin{itemize}
\item[(i)] the area of\/ $\mathbb{S}^2$ with
respect to the induced metric
$\Area(\mathbb{S}^2,f^*g_{\mathbb{S}^n})$
is an integer multiple of $4\pi;$
\item[(ii)] $n$ is even, $n=2m,$ and
\item[(iii)]
$$
\Area(\mathbb{S}^2,f^*g_{\mathbb{S}^n})\geqslant2\pi m(m+1).
$$
\end{itemize}
\end{Theorem}

\begin{Definition}\label{degree-definition}
If $\Area(\mathbb{S}^2,f^*g_{\mathbb{S}^n})=4\pi d,$ then we say that $f$
is of {\em harmonic degree}\/ $d.$
\end{Definition}

Theorem~\ref{Calabi-Barbosa} implies
the following proposition.

\begin{Proposition}\label{sufficient-maps}
Extremal metrics (possibly with conical singularities)
on $\mathbb{S}^2$ are induced by linearly full
harmonic immersions with branch points 
$f:\mathbb{S}^2\longrightarrow\mathbb{S}^{2m}.$
The harmonic degree $d$ of such an immersion satisfies
the inequality $d\geqslant\frac{m(m+1)}{2}.$ 
\end{Proposition}

Indeed, any non-linearly full immersion can be a viewed as 
linearly full immersion to a corresponding subsphere.

As we know from 
Theorem~\ref{metrics-harmonic-maps-s2-branch},
a metric $g$ induced by a linearly full
harmonic immersion with branch points 
$f:\mathbb{S}^2\longrightarrow\mathbb{S}^{2m}$
is extremal for the eigenvalue 
$\lambda_{N(2)}$, where
$$
N(2):=N(2;\,\mathbb{S}^2,g).
$$
To simplify notation, we will use this convention in the remainder of the paper.

A priori, it is not clear how to find 
$N(2)$ for a metric on a surface induced by
a harmonic immersion into a sphere. 
The following lower bound on $N(2)$  due to Ejiri is a 
particular case of \cite[Theorem C]{Ejiri1998} and is crucial 
for our considerations.

\begin{Theorem}[Ejiri~\cite{Ejiri1998}]\label{Ejiri}
Let $f:\mathbb{S}^2\longrightarrow\mathbb{S}^{2m}$ 
be a linearly full harmonic map
of harmonic degree $d>1$ of $\mathbb{S}^2$
to the standard unit sphere $\mathbb{S}^{2m}.$
Then
$$
N(2)\geqslant d+1.
$$
\end{Theorem}
We note that while the condition $d>1$ is not explicitly stated in the formulation
of Theorem C in \cite{Ejiri1998}, it is assumed and used in its proof. 
For details, see Appendix \ref{sec:Ejiri}.

\section{Proof of Theorem~\ref{maintheorem}}\label{mainproof}
As was mentioned in the introduction, it is known
that 
$
\Lambda_k(\mathbb{S}^2)=8\pi k
$
for $k=1$  (\cite{Hersch1970}),
$k=2$ (\cite{Nadirashvili2002, Petrides2014})
and $k=3$ (\cite{Nadirashvili-Sire2015q2}).

Let us first show that  for $k=1$  a round metric is, up to rescaling, 
the unique maximal metric on the sphere among all metrics that are smooth 
outside a finite number of  conical singularities.
This is well-known for smooth metrics. Indeed, it was proved 
in~\cite{Montiel-Ros1986} (see also \cite{ElSoufi1986})
that there exists at most one metric in a conformal class  
admitting an isometric immersion into a unit sphere by the 
first eigenfunctions.  Given that there is only one conformal 
structure on  $\mathbb{S}^2$, the uniqueness of the maximal metric 
readily follows from Theorem~\ref{Nadirashvili-ElSoufi-Ilias}.
It remains to show that a maximal metric for the first eigenvalue 
on the sphere is smooth. It follows from Proposition~\ref{sufficient-maps} 
that a maximal metric is induced by a linearly full harmonic immersion 
with branch points 
$f:\mathbb{S}^2\longrightarrow\mathbb{S}^{2m}$
of harmonic degree $d.$ Since it is maximal for the
first eigenvalue, we have $N(2)=1.$ It follows from Theorem~\ref{Ejiri}
that if $d>1$ then $1\geqslant d+1,$ which  is impossible.
Therefore,  $d=1$,  and hence $m=1$  by Proposition~\ref{sufficient-maps}. 
Moreover, $\Area(\mathbb{S}^2,f^*g_{\mathbb{S}^2})= 4\pi$, which implies 
that $f$ is a one-sheeted covering. In particular, it  does not have branch 
points, and thus  a  maximal metric is smooth.
 
We now  prove \eqref{main:isop}  by induction in  $k.$  Indeed, the result
holds for $k=1$ as was mentioned  above. Suppose that it holds for $j=1,...,k-1$.
Assume, by contradiction,  that 
$\Lambda_k(\mathbb{S}^2)>8\pi k$. Then by Theorem \ref{thm:petrides} 
there exists a maximizing metric $g$ which is smooth possibly with the 
exception of a finite number of conical singularities such 
that $\Lambda_k(\mathbb{S}^2)=\bar \lambda_k(\mathbb{S}^2, g)$.
By Proposition \ref{sufficient-maps} the metric $g$ 
is  induced on $\mathbb{S}^2$
by a linearly full harmonic immersion (possibly with branch points)
$f: \mathbb{S}^2\longrightarrow\mathbb{S}^{2m}$
of harmonic degree $d\geqslant 1$.
Note that by Theorem \ref{metrics-harmonic-maps-s2-branch}, 
\begin{equation}
\label{xxx}
\lambda_k(\mathbb{S}^2,g)=2.
\end{equation}
At the same time,
\begin{equation}
\label{triv}
\lambda_k(\mathbb{S}^2,g)=2>\lambda_j(\mathbb{S}^2,g)
\end{equation}
for all $j=0,1,\dots,k-1$, and therefore 
$$N(2)=k.$$ Indeed, it follows from  the induction hypothesis that for any $j<k$,
$$\lambda_k(\mathbb{S}^2,g) \Area(\mathbb{S}^2,g) = %
\Lambda_k(\mathbb{S}^2) > 8\pi k > 8\pi j = %
\Lambda_j(\mathbb{S}^2)\geqslant \lambda_j(\mathbb{S}^2,g) \Area(\mathbb{S}^2,g),$$
which immediately implies inequality \eqref{triv}.

Consider  the case $d>1.$
Then Theorem~\ref{Ejiri} implies that 
\begin{equation}\label{N-d}
N(2)=k\geqslant d+1.
\end{equation}
On the other hand, by Definition~\ref{degree-definition} 
we have
$$
\Area(\mathbb{S}^2,g)=4\pi d.
$$
Therefore, combining \eqref{xxx}  and \eqref{N-d} we get
$$
\Lambda_k(\mathbb{S}^2)=\bar{\lambda}_{k}(\mathbb{S}^2, g)=%
2\Area(\mathbb{S}^2, g)=8\pi d<8\pi k,
$$
which contradicts the assumption 
$\Lambda_k(\mathbb{S}^2)>8\pi k$.

Finally, consider the case $d=1.$ Then 
$\Area(\mathbb{S}^2, g)=4\pi$ and
$$
\bar{\lambda}_k(\mathbb{S}^2, g)=8\pi.
$$
However, in this case the  metric $g$ could be maximal only for $k=1,$
but we have assumed that  $k>1.$ 
Therefore $\Lambda_k(\mathbb{S}^2)=8\pi k$, for any $k\geqslant 2$, and 
it follows from Remark \ref{gluing}  that the equality is attained in the  
limit by a sequence of metrics degenerating to a union of $k$ touching 
identical round spheres. 
This completes the proof of the theorem. $\Box$ 

\appendix

\section{Theorem \ref{Ejiri}: outline of the proof}
\label{sec:Ejiri}

The goal of this Appendix is to present the main ideas of the proof 
of Theorem \ref{Ejiri} of~\cite{Ejiri1998} which is a key ingredient in  
the proof of Theorem \ref{maintheorem}. The techniques used 
in~\cite{Ejiri1998} are not well-known in the spectral geometry community, 
but they appear to be quite powerful (see \cite{Nayatani-Shoda2017q} for 
another recent application), and we believe it is worth discussing them in 
some detail. In particular, for the convenience of the reader, we recall  
a number of notions from algebraic geometry that are used in the proof of 
Theorem~\ref{Ejiri}.

\subsubsection*{I. Isotropy} Recall that there exists a unique complex 
structure on the sphere $\mathbb{S}^2$.
Let $z$ be a local complex coordinate, and denote by $\dz$ and $\dzb$ 
derivatives with respect to $z$ and $\bar z$ respectively. One can easily 
check that the definitions below do not depend on the choice of $z$. For 
two vectors $u,v\in \mathbb{C}^n$ we use $(u,v) = \sum\limits_{i=1}^nu_iv_i$ 
to denote $\mathbb{C}$-bilinear inner product and use 
$\langle u,v\rangle = (u,\bar v) = \sum\limits_{i=1}^nu_i\bar{v}_i$ 
to denote the Hermitian inner product. The symbol $\perp$ refers to 
orthogonality in the Hermitian metric. 

\begin{Definition}
A subspace $V$ of $\mathbb{C}^n$ is called {\em isotropic} if $(v_1,v_2) = 0$ 
for all $v_1,v_2\in V$. Equivalently, $V$ is isotropic if $V\perp\bar V$ 
with respect to the Hermitian inner product on $\mathbb{C}^n$.
\end{Definition}

\begin{Definition}
\label{isotropic}
A branched conformal immersion 
$x\colon \mathbb{S}^2 \to \mathbb{S}^n\subset \mathbb{R}^{n+1}$ 
is called {\em isotropic} if for all $k+l\geqslant 1$ one has 
$(\dz^k x, \dz^l x) = 0$, i.e. for any point $p\in \mathbb{S}^2$, the  
complex derivatives of $x$ at $p$ span an isotropic subspace in 
$\mathbb{C}^{n+1}$.
\end{Definition}
Calabi~\cite{Calabi1967} proved the following theorem.
\begin{Theorem}
Any linearly full branched minimal immersion 
$x\colon\mathbb{S}^2\to\mathbb{S}^n$ is isotropic. 
\end{Theorem}
\begin{Remark}
Note that in the literature, isotropic minimal immersions to $\mathbb{S}^{4}$ 
are sometimes referred to as {\em superminimal}. The notion of 
superminimal immersion was introduced in the paper of Bryant~\cite{Bryant} 
though the definition there differs from Definition~\ref{isotropic}.
\end{Remark}

As a corollary one obtains the statement (ii) of 
Theorem~\ref{Calabi-Barbosa}, i.e. $n =2m$.
\begin{proof}
Fix a point $p\in \mathbb{S}^2$. Let $G$ be an isotropic subspace spanned 
by complex derivatives of $x$ at $p$. Then $G$ is isotropic, i.e. $G\perp \bar G$.
Moreover, as $x$ sends $\mathbb{S}^2$ to the unit sphere, one 
has $(x,x)=1$. Applying $\dz^k$ and $\dzb^k$ to this relation and 
using the isotropy of $G$ and $\bar G$ one obtains $x(p)\perp G$ and 
$x(p)\perp \bar G$. 

Suppose that there exists a constant vector $u\in \mathbb{R}^{n+1}$ such 
that $u\perp \mathbb{C} x(p)\oplus G\oplus \bar G$. Then $(x,u)$ has a 
zero of infinite order at $p$ and by the unique continuation theorem 
$(x,u)\equiv 0$,  which contradicts the fact that $x$ is linearly full. 
Therefore, $\mathbb{C}^{n+1} = \mathbb{C} x\oplus G\oplus \bar G$, i.e. 
$n=2m$.
\end{proof}

\subsubsection*{II. Directrix curve}  To any linearly full harmonic 
map $x\colon \mathbb{S}^2 \to \mathbb{S}^{2m}$, one can canonically 
associate a holomorphic curve 
$\Phi\colon \mathbb{S}^2\to\mathbb{CP}^{2m}$  called the directrix curve. In 
order to introduce  it we need the two auxiliary  definitions below.

\begin{Definition}
For a holomorphic curve $\Psi\colon \mathbb{S}^2 \to\mathbb{CP}^{n}$,  
the {\em $k$-th associated curve} $\Psi_k$ is the holomorphic map 
$\Psi_k\colon \mathbb{S}^2 \to Gr_{k+1,n+1}(\mathbb{C})$ to a 
complex Grassmannian,  which  in terms of any local holomorphic 
lift $\psi\colon U \to \mathbb{C}^{n+1}$ is given by a $k$-subspace 
spanned by $\psi,\dz\psi,\ldots,\dz^{k}\psi$.
\end{Definition}

\begin{Definition}
A linearly full holomorphic curve 
$\Psi\colon \mathbb{S}^2 \to\mathbb{CP}^{2m}$ is called {\em isotropic} 
if its $(m-1)$-th associated curve takes values in isotropic $m$-planes 
in $\mathbb{C}^{2m+1}$.
\end{Definition}

For each isotropic curve $\Psi\colon \mathbb{S}^2\to\mathbb{CP}^{2m}$,  
one can construct a pair of linearly full harmonic maps 
$\pm x_{\Psi}\colon \mathbb{S}^2\to\mathbb{S}^{2m}$ in the 
following way. In view of isotropy,  
$\Psi_{m-1}\perp\bar\Psi_{m-1}$ as $m$-dimensional subspaces of 
$\mathbb{C}^{2m+1}$. Thus, $\Psi_{m-1}\oplus\bar\Psi_{m-1}$ spans a 
$2m$-dimensional subspace, which is invariant under conjugation. 
Therefore, its orthogonal complement is also invariant under conjugation, 
and we define $\tilde x_{\Psi}\colon \mathbb{S}^2\to\mathbb{RP}^{2m}$ as 
a unique real line in $(\Psi_{m-1}\oplus\bar\Psi_{m-1})^\perp$.  
Barbosa proved~\cite{Barbosa1975} that such $\tilde x_\Psi$ is harmonic.
As $\mathbb{S}^2$ is simply connected, $\tilde x_{\Psi}$ can be lifted to 
a pair of linearly full harmonic maps 
$\pm x_{\Psi}\colon \mathbb{S}^2\to\mathbb{S}^{2m}$.

\begin{Proposition} [Barbosa \cite{Barbosa1975}]
For any linearly full harmonic map 
$x\colon \mathbb{S}^2 \to \mathbb{S}^{2m}$ there exists unique 
isotropic curve $\Phi\colon \mathbb{S}^2\to\mathbb{CP}^{2m}$ such that 
$\pm x_\Phi = x$.
\end{Proposition}

The curve $\Phi$ is called {\em directrix curve} of the linearly full 
harmonic map $x$.
In what follows we do not distinguish between $x$ and $-x$, 
see Remark~\ref{Cptobvious}.

\begin{proof}[Idea of the proof.]
For each $i =1,\ldots, m$,  define $G_i$ inductively as follows: 
$G_1 = \dzb x$ and 
$$
G_i = \dzb^i x + \sum\limits_{j=1}^{i-1} a_i^j G_j,
$$
where the coefficients $a_i^j$ are chosen in such a way that 
$\langle G_i, G_j\rangle = 0$ for $i\ne j$. Then $\Phi = \pi\circ G_m$, 
where $\pi\colon \mathbb{C}^{2m+1}\backslash\{0\}\to \mathbb{CP}^{2m}$ 
is the natural projection.
\end{proof}

\subsubsection*{III. Harmonic degree} In this section we give another 
interpretation of the harmonic degree of an isotropic  branched minimal immersion, 
see  Definition \ref{degree-definition}.

Let $\xi_{k,n}\colon Gr_{k,n}(\mathbb{C})\to \mathbb{CP}^{\binom{n}{k}-1}$ 
denote the Pl\"ucker embedding, see e.g.~\cite{GH}. Let $\tilde d$ be the 
algebraic degree of the curve 
$$
\xi_{m,2m+1}\circ\Phi_{m-1}\colon  \mathbb{S}^2 \to 
\mathbb{CP}^{\binom{2m+1}{m}-1}.
$$
Then one has the following theorem.
\begin{Theorem}[Barbosa~\cite{Barbosa1975}]
The degree $\tilde d$ is even. Moreover, it is equal to twice the harmonic 
degree $d$. 
\end{Theorem}

\subsubsection*{IV. Action of $SO(2m+1,\mathbb{C})$} 
Barbosa~\cite{Barbosa1975} noticed that the group $SO(2m+1,\mathbb{C})$ 
naturally acts on the space of linearly full isotropic branched minimal immersions. 
Indeed, let $A\in SO(2m+1,\mathbb{C})$ and let $\Phi$ be the directrix curve 
of $x\colon  \mathbb{S}^2\to\mathbb{S}^{2m}$ (see subsection II). Then $A\Phi$ 
is an isotropic  curve in $\mathbb{CP}^{2m}$, and therefore it is the 
directrix curve of some  linearly full isotropic branched 
minimal immersion that we denote 
by $Ax$. 

In the same way, $SO(2m+1,\mathbb{C})$ acts on $m$-dimensional subspaces,  
and it is easy to see that $(A\Phi)_{m-1} = A(\Phi_{m-1})$. Thus, this 
action preserves the harmonic degree of the immersion $x$. This fact will 
be used  in Proposition~\ref{Ejiri:prop}.

\subsubsection*{V. Conformally covariant Schr\"odinger operators} 
To each conformal map $x\colon  \mathbb{S}^2 \to\mathbb{S}^{2m}$ one can 
associate~\cite{MR} a conformally covariant Schr\"odinger operator $L_x$ in 
the following way. For any metric $g$ on $\mathbb{S}^2$ we set
$$
L_{x,g} = \Delta_g - |\nabla x|_g^2,
$$
where $|\nabla x|_g^2 = \sum_{i=1}^{2m+1}|\nabla x_{i}|_g^2$. It is easy to 
see that for any function $\omega\in C^\infty( \mathbb{S}^2)$ one 
has $L_{x,e^{2\omega}g} = e^{-2\omega}L_{x,g}$. Thus, the corresponding 
quadratic form
$$
Q_x(u) = \int_{\mathbb{S}^2}|\nabla u|_g^2 - |\nabla x|_g^2 u^2\,dv_g
$$
is conformally invariant. This fact motivates the following definition.

\begin{Definition}
The {\em index} $\ind(x)$ of a conformal map 
$x\colon  \mathbb{S}^2\to\mathbb{S}^{2m}$ is defined to be the index 
of $L_{x,g}$,~i.e. the number of negative eigenvalues of~$L_{x,g}$. 

Similarly, the {\em nullity} $\nul(x)$ of a conformal map 
$x\colon  \mathbb{S}^2\to\mathbb{S}^{2m}$ is defined as the nullity of 
$L_{x,g}$,~i.e. $\dim\ker L_{x,g}$. Both index and nullity do not depend 
on the choice of the metric $g$. 

\end{Definition}
\begin{Remark}
\label{Cptobvious}
Clearly, $\ind(x) = \ind(-x)$ and $\nul(x) = \nul(-x)$.
\end{Remark}
\begin{Remark}
\label{remark:index}
An elementary computation shows that for $g = x^*g_{\mathbb{S}^{2m}}$ 
the operator $L_{x,g}$ takes the form
$$
L_{x,g} = \Delta_g - 2,
$$
i.e. $\ind(x)$ coincides with the number $N(2)$ appearing in the statement of 
Theorem \ref{Ejiri}.
\end{Remark}

\subsubsection*{VI. Proof of Theorem~\ref{Ejiri}} 
The following proposition summarizes some of the results of  the 
paper~\cite{Ejiri1998}.
\begin{Proposition}[Ejiri~\cite{Ejiri1998}]
\label{Ejiri:prop}
Let $x\colon \mathbb{S}^2\to\mathbb{S}^{2m}$, $m>1$,  be a linearly full 
isotropic branched minimal immersion of harmonic degree $d>1$. Then there exists 
a one-parameter subgroup $A(t)\subset SO(2m+1,\mathbb{C})$ such that
\begin{itemize}
\item[1)] As $t\to\infty$, the one-parameter family $x_t =A(t)x$ 
converges uniformly to a linearly full isotropic branched minimal immersion 
$x_\infty\colon \mathbb{S}^2\to\mathbb{S}^{2m-2}$ of harmonic degree 
$d$, where $\mathbb{S}^{2m-2}\subset \mathbb{S}^{2m}$ is an equatorial sphere;
\item[2)] $\nul(x) = \nul(x_t) \leqslant \nul(x_\infty)$;
\item[3)] $\ind(x) = \ind(x_t)\geqslant \ind(x_\infty)$.
\end{itemize}
\end{Proposition}
In the first statement we are using the degree-preserving property of the 
action of $SO(2m+1, \mathbb{C})$ which was explained in subsection IV. Note 
also that the third statement follows from the second  one by continuity 
of eigenvalues of the Schr\"odinger operator $L_{x,g}$ with respect to 
continuous  deformations of the potential.

If $m>2$, one can apply the previous proposition to $x_\infty$ and obtain a 
new  isotropic branched minimal immersion to $\mathbb{S}^{2m-4}$. Continuing this 
process, we arrive at a linearly full isotropic branched minimal immersion 
$y\colon \mathbb{S}^2\to \mathbb{S}^2$ of harmonic degree $d$.  
Moreover, $\ind(x) \geqslant \ind(y)$.  At the same time,  the 
following result holds:

\begin{Theorem}[Montiel and Ros \protect{\cite[Theorem 21]{MR}}; 
Nayatani~\cite{Nayatani}]
Let $y\colon \mathbb{S}^2\to\mathbb{S}^2$ be a conformal map of degree $d>1$. 
Then $\ind(y)\geqslant d+1$.
\end{Theorem}
Therefore, given a linearly full harmonic map 
$x\colon  \mathbb{S}^2 \to\mathbb{S}^{2m}$ of degree $d>1$, in view of 
Remark \ref{remark:index} we have 
$N(2)=\ind(x)\geqslant \ind(y)\geqslant d+1$, which completes the proof of 
Theorem \ref{Ejiri}. \qed

\end{document}